\documentclass[10pt,twoside]{article}

\usepackage{amsmath}
\usepackage{amsfonts}
\usepackage{amssymb}
\usepackage{amscd}
\usepackage{amsthm}
\usepackage{amsbsy}
\usepackage{graphicx}
\usepackage{bm}

\def\@oddhead{\hfill \shorttitle \hfill \thepage}
\def\@evenhead{\thepage \hfill \shortauthor \hfill}
\def\@oddfoot{}
\def\@evenfoot{}

\newtheorem{thm}{Theorem}[section]

\newtheorem{lem}[thm]{Lemma}
\newtheorem{prop}[thm]{Proposition}

\newtheorem{defn}[thm]{Definition}
\theoremstyle{remark}

\numberwithin{equation}{section}

\newcommand{\R}{\mathbb R}

\def\XXint#1#2#3{{\setbox0=\hbox{$#1{#2#3}{\int}$}
     \vcenter{\hbox{$#2#3$}}\kern-.5\wd0}}

\newcommand{\HH}{\mathbb H}

\newcommand{\RR}{\mathbb R}
\newcommand{\ZZ}{\mathbb Z}
\newcommand{\del}{\partial}

\newcommand{\Ric}{\mathrm{Ric}}

\newcommand{\calC}{{\mathcal C}}
\newcommand{\calD}{{\mathcal D}}

\newcommand{\calO}{{\mathcal O}}

\newcounter{mnotecount}[section]
\let\oldmarginpar\marginpar
\renewcommand\marginpar[1]{\-\oldmarginpar[\raggedleft\footnotesize #1]%
{\raggedright\footnotesize #1}}

\date{}
\title{\ \\[0.4cm] \ \\ \bf  Ricci flow in two dimensions }
\author{James Isenberg \footnote{University of Oregon, Department of Mathematics, Eugene, OR 97403}\hspace{2mm}
Rafe Mazzeo \footnote{Stanford University, Department of Mathematics, Stanford, CA 94305} 
\hspace{2mm} and Natasa Sesum \footnote{Rutgers University, Department of Mathematics, New Brunswick, NJ 08854}}
\begin{document}

\maketitle


\thispagestyle{empty}

\begin{abstract}
\vskip 3mm\footnotesize{

\vskip 4.5mm
\noindent
Ricci flow on two dimensional surfaces is far simpler than in the higher dimensional cases. This presents
an opportunity to obtain much more detailed and comprehensive results. We review the basic facts about this 
flow, including the original results by Hamilton and Chow concerning Ricci flow on compact surfaces.
The rationale for this paper, however, is especially to survey recent work concerning this flow on open surfaces, including various 
classes of both complete and incomplete surfaces, where a number of striking new phenomena have been observed.

\vspace*{2mm}
\noindent{\bf 2000 Mathematics Subject Classification: 53C44, 30F45}

\vspace*{2mm}
\noindent{\bf Keywords and Phrases: Ricci flow, uniformization, complete surfaces}}

\end{abstract}
\vspace*{-8.2cm}
\vspace*{2mm} \noindent\hspace*{82mm}
\begin{picture}(41,10)(0,0)\thicklines\setlength{\unitlength}{1mm}
\put(0,2){\line(1,0){41}} \put(0,16){\line(1,0){41}}%
\put(0,12.){\sl \copyright\hspace{1mm}Higher Education Press}
\put(0,7.8){\sl \hspace*{4.8mm}and International Press}%
\put(0,3.6){{\sl \hspace*{4.8mm}Beijing-Boston} }
\end{picture}

\vspace*{-17.6mm}\noindent{{\sl The title of\\
This book*****}\\ALM\,?, pp.\,1--?} \vskip8mm
\vspace*{7.6cm}

\section{Introduction}
Ricci flow is by now so well-established as a fundamental tool in geometric analysis that it barely needs an introduction. 
Beyond its spectacular applications to topology, its development has led to many innovations in the analysis of nonlinear 
parabolic equations, and has illustrated how geometric insight can provide significant new tools for their study.  There
are now several excellent introductions to the Ricci flow in higher dimensions 
\cite{Brendle-book, CK, CCGGGIIKLLN1, CCGGGIIKLLN2, MT, Tbook}.  Our focus in this survey, however, is the lowest dimensional 
case where Ricci flow makes sense, namely for surfaces. The expectation that the mathematics of two dimensional Ricci 
flow should  be especially interesting and useful is suggested by the analogous situation for extrinsic curvature flows: 
there are particularly complete results for curve-shortening flow, and these elucidate and motivate many aspects of higher 
dimensional mean curvature flow. 

In two dimensions, the Ricci curvature for a metric $g$ is equal to $\frac12 R g$, where $R$ is the scalar curvature (or
twice the Gauss curvature). Hence the Ricci flow equation for surfaces takes the form
\begin{equation}
\del_t g(t) = (\rho - R(t)) g(t),
\label{eq:rfs0}
\end{equation}
where $\rho$ is any normalizing constant. The choice of $\rho$ is not too important, and indeed there is a simple transformation which 
interchanges solutions for one value of $\rho$ with solutions for any other value of $\rho$; however, as we discuss below, 
in certain situations there are particular choices of $\rho$ which simplify the analysis somewhat. 

From equation \eqref{eq:rfs0}, it is clear that two dimensional Ricci flow solutions, if they exist, remain within a 
conformal class, and it is also clear that this flow
is the same as the Yamabe flow on surfaces. This has the very important consequence that two dimensional Ricci 
flow reduces to a scalar problem. Consequently, the maximum principle can be used in very direct ways and leads
to particularly strong conclusions.  Since Riemann surfaces are automatically K\"ahler, one can also consider two dimensional 
Ricci flow as an instance of the K\"ahler-Ricci flow, involving a (different) scalar parabolic equation. This alternate perspective 
leads to further insights and to new techniques which can be applied to this problem.

The most celebrated application of Ricci flow in higher dimensions is to the diffeomorphic classification of 
manifolds.  For surfaces, however, this is not so interesting since the classification of surfaces, both closed
and open, is well-known and follows from the general uniformization theorem, which states that any Riemannian surface 
$(M^2,g)$ (closed or open, complete or incomplete) admits a conformally related metric which is complete 
and has constant curvature, and further this `uniformizing metric' is unique up to scaling except in the special case 
$M^2 = S^2$. The classical proof of uniformization uses complex function theory, but we note that it can now be proved 
by fairly simple PDE methods \cite{Maz-Tay}.

Based on the uniformization theorem, an optimist might conjecture that if $(M^2,g_0)$ is any (not necessarily closed) Riemannian 
surface, then the Ricci flow starting at $g_0$ exists for all time and converges to the uniformizing constant curvature metric in 
the conformal class of $g_0$ on $M^2$, or some multiple of it.  This is most plausible if $M^2$ is compact or if
$(M^2,g_0)$ is complete, but if $(M^2,g_0)$ is incomplete the statement must be modified since as we discuss in the last
section of this paper, the flow is not uniquely determined then.  Our goal here is to discuss how much of this is known
to be true. For compact surfaces, the theory is complete and Ricci flow provides a new proof of the uniformization theorem, 
albeit with some twists which complicate the proof in the spherical setting. Much is also known in the complete case.
There are also many fascinating issues related to setting up this flow in the incomplete setting, and we discuss some 
of these below as well.  

In any case, the geometric interest in studying this set of problems is not so much one of proving the existence of canonical 
metrics, but rather of proving the stability of uniformizing metrics. More specifically, if $\bar{g}$ is a complete metric of 
constant curvature, then we seek to understand the size of the neighbourhood of $\bar{g}$ in the space of all metrics on 
$M^2$ on which the flow has $\bar{g}$ as an attractor basin. In other words, we seek to understand whether there are simple
quantitative geometric conditions on a metric $g_0$ so that the Ricci flow with initial metric  $g_0$ converges to a constant 
multiple of $\bar{g}$. 

To conclude this introduction, all three authors wish to record our debt of gratitude to Rick Schoen for
his inspired and tireless leadership in geometric analysis. We each have profited greatly from our interactions
with him over the years, and are very pleased to be able to contribute to this volume of papers in his honor.

The authors are supported by the NSF grants PHY-0652903,  DMS-0805529 and DMS-0905749, respectively. 

\section{General considerations}
\label{GenCon}
As noted above, since in two dimensions $\Ric(g) = \frac12 R \, g$ where $R$ is the scalar curvature, it follows that on 
surfaces the general Ricci flow equation with volume-normalizing term $\del_g g(t) = -2\Ric(g(t)) + \rho g(t)$ 
takes the particularly simple form \eqref{eq:rfs0}. This flow clearly preserves the conformal class of
the metric, and hence the evolution equation can be recast as a scalar parabolic equation.  To do this, we recall the transformation law
\begin{equation}
\Delta_{0} \phi - \frac12 R_{0} + \frac12 R e^{2\phi} = 0
\label{eq:trsc}
\end{equation}
relating the scalar curvatures $R_0$ and $R$ of two conformally related metrics $g = e^{2\phi}g_0$. Then if we write 
$g(t) = u(t)g_0$,  we see that \eqref{eq:rfs0} is equivalent to the scalar equation (of ``fast diffusion'' type) 
\begin{equation}
\del_t u = \Delta_{g_0} \log u + \rho u - R_0, \qquad u(0) = u_0,
\label{eq:rfs}
\end{equation}
written here as an initial value problem.
Already we see a simplification which results from working in this low dimension: in all higher dimensions the Ricci flow equation is not 
parabolic because of its invariance under the diffeomorphism group; however  in this case it is and consequently in proving well-posedness no gauge conditions need to be imposed

Combining \eqref{eq:trsc} and \eqref{eq:rfs}, and using  $\phi = \frac12 \log u$, we obtain the simple and
useful relationship
\begin{equation}
\del_t u = (\rho-R) u \quad \Longleftrightarrow \quad \del_t \log u = \rho - R.
\label{eq:derivu}
\end{equation}
This equation leads  us to focus attention on controlling the function $R$ as $t$ increases, since it follows from  \eqref{eq:derivu} 
that appropriate bounds and convergence of this function directly imply analogous bounds and convergence for $u$ itself. To 
study $R$ directly, it is useful to work with its evolution equation, which one derives by differentiating \eqref{eq:trsc} and 
presuming that $g(t)$ is a solution of \eqref{eq:rfs0} on some interval $0 \leq t < T$; one obtains
\begin{equation}
\label{eq-R1}
\del_t R = \Delta_{g(t)} R + R(R-\rho).
\end{equation}
An immediate important consequence of equation \eqref{eq-R1} is that the minimum of $R$ is non decreasing  in time in any situation 
in which  the minimum principle can be applied. For example, suppose that $M^2$ is compact and that $\rho$
is chosen to be the average of $R$ over $M^2$ (which is independent of $t$ by Gauss-Bonnet, though the $t$-independence
can be deduced easily from \eqref{eq-R1} and \eqref{eq:evolaf} below). Then the function  $R_{\min}(t) = \inf_{M} R(z,t)$, which 
depends only on $t$,  satisfies
\[
\frac{d}{dt} R_{\min} \geq R_{\min} (R_{\min} - \rho).
\]
Note that $R_{min}(t)$ is only Lipschitz, but the precise meaning of this differential inequality is that the term on the left is 
the lim inf of the forward difference quotient of $R_{\min}(t)$. 
Since by definition $R_{\min} \leq \rho$, one sees immediately that if $\rho \le 0$, then the right hand side 
is nonnegative, and the claim holds. If instead $\rho > 0$, then a similar argument applied to the difference
$R_{\min} - r(t)$, where $r(t)$ is chosen so that $\frac{d}{dt} r(t) = r(t) (r(t) - \rho)$ and $r(0) = R_{\min}(0)$, leads to  the same conclusion.  
In many noncompact settings, it is possible to apply similar reasoning. Interestingly, in all situations,
obtaining an upper bound for $R$ is more difficult.

We conclude this section with a remark about the choice of the constant $\rho$ in \eqref{eq:rfs0}. For all
considerations which are local in $t$, the choice of $\rho$ is irrelevant since it can be altered to any other
value with a simple `gauge' transformation.  For long-time considerations, however, it does have an influence.
For example, if $dA$ is the area form for $g(t)$, then 
\begin{equation}
\del_t (dA) = (\rho-R) dA.
\label{eq:evolaf}
\end{equation}
Using this, if $M^2$ is compact and if $\mathcal{A}(t) := \int_{M^2} dA$ is the area at time $t$, then
\[
\frac{d}{dt} \mathcal{A}(t) = \int_{M^2} (\rho - R)\, dA = \rho \mathcal{A}(t) - 4\pi \chi(M^2).
\]
Hence if $\rho = 4\pi \chi(M^2)/\mathcal{A}(0)$, then $\mathcal{A}(t) \equiv \mathcal{A}(0)$;  i.e., the area is independent of $t$.
Other choices of $\rho$ lead to finite time collapse.  In certain situations in which $M^2$ is non compact, we 
can usefully make sense of this same calculation.  

\section{Compact surfaces}
\label{compact}
The earliest results concerning the long-time existence and convergence properties for Ricci flow solutions 
in two dimensions have involved geometries on compact manifolds. As already noted,  for compact surfaces it is reasonable 
to choose $\rho = 4\pi \chi(M)/\mathcal{A}(0)$, so that the area remains constant in $t$. The complete story for solutions of \eqref{eq:rfs0} 
with this normalization on compact surfaces is presented in a sequence of two papers written during the late 1980's,  by 
Hamilton \cite{Ha1} and Chow \cite{Ch1}, respectively, which contain the following two results:
\begin{thm}[Hamilton \cite{Ha1}]
\label{thm-ham}
Let $(M^2,g_0)$ be compact. If $\rho \le 0$, or if $R(0) \geq 0$ on all of $M^2$, then the solution to (\ref{eq:rfs0}) exists 
for all $t \geq 0$ and converges to a metric of constant curvature.
\end{thm}
\begin{thm}[Chow \cite{Ch1}]
\label{thm-chow}
If $g_0$ is any metric on $S^2$, then its evolution under (\ref{eq:rfs0}) develops positive scalar 
curvature in finite time, and hence by Theorem \eqref{thm-ham} converges to the round metric as $t \nearrow \infty$. 
\end{thm}
In the remainder of this section we discuss various features of the proofs of these theorems, and refer to these two
papers for complete details.

\subsection{Long time existence of the flow}
As with most evolution equations, one proves that solutions exist for all $t \geq 0$ by combining a short-time
existence (and uniqueness) result with a priori bounds which show that solutions cannot develop singularities
in finite time.  For compact surfaces, short-time existence is proved in the standard way for the scalar parabolic 
equation \eqref{eq:rfs}: by setting it up as a fixed point problem for a contraction mapping. This mapping
is obtained by applying the fundamental solution of the linearization at any given $u_0 > 0$ to \eqref{eq:rfs};
it is a contraction on any sufficiently short time interval. 

To obtain the needed a priori estimates, it suffices to find upper and lower bounds for $R(t)$, because using such 
bounds, together with \eqref{eq:derivu}, we determine that $u$ itself is bounded on any finite time-interval $[0,T)$.
It then follows from  standard parabolic Schauder estimates that $u$ is bounded in
 H\"older spaces with arbitrarily
high norms, whence we obtain compactness by the Arzela-Ascoli theorem. This is the simplest occurrence 
of a standard theme in the study of the higher dimensional flows, that the flow develops singularities only 
when the curvature blows up.  

The minimum principle certainly holds here, so the argument we have sketched above (in \S \ref{GenCon}) to show that $R_{\min}$ is
nondecreasing applies.  The upper bound for $R$ is less straightforward, and one of Hamilton's innovations 
was to show that this supremum can be controlled using a so-called potential function. This is, by definition, 
the unique solution $f$ to the equation
\[
\Delta f = R - \rho,
\]
normalized to have integral $0$. (Note that this equation always admits a solution since $\int (R - \rho) = 0$).  

A computation shows that $f$ satisfies the evolution equation
\[
\del_t f = \Delta f + rf - b,
\]
where $b = b(t) := \frac{1}{\mathcal{A}\,(M))}\int_M |\nabla f|^2\, dA$ depends only on $t$. Hamilton's key observation
is that the function $h := \Delta f + |\nabla f|^2$ satisfies an evolution equation
\[
\del_t h = \Delta h - 2|M_f|^2 + \rho h,
\]
where $M_f = \nabla^2 f - \frac12 \Delta f g$ is the trace-free (covariant) Hessian of $f$.  Using the maximum principle to estimate $h$, 
and the relationship $R = h - |\nabla f|^2 + \rho$, we obtain the  bound
\begin{equation}
\label{RBound}
-C_1 \le R \le C_2 e^{\rho t} + \rho
\end{equation}
for some constants $C_1, C_2$, which gives uniform bounds on any finite time interval.  This yields a priori bounds 
for $u$ on any $[0,T)$, and hence proves long-time existence of the flow. 

\subsection{Convergence} 
In verifying that Ricci flow solutions (with the area-preserving normalization given above) on compact surfaces always converge to a constant curvature metric, there are four separate cases to consider:
$\rho < 0$, $\rho = 0$, $\rho > 0$ with $R(0)$ everywhere nonnegative, and
finally $\rho > 0$ with $R(0)$ changing sign. The first two are quite simple, the third requires some significant new ideas,
and the fourth, which is Chow's contribution, involves further refinements. 

\medskip
\noindent{$\bf{\rho < 0}:$ }
In this case, it follows from the arguments in \S \ref{GenCon} that $R_{\min}(t)$ is strictly increasing and converges to $\rho$ from 
below. On the other hand, by \eqref{RBound}, $R_{\max}(t)$ decreases exponentially to $\rho$ from above, i.e.
\[
\rho - C e^{\rho t} \leq R(t) \leq \rho + C e^{\rho t}
\]
for all $t \geq 0$; hence $R \to \rho$ at an exponential rate. Integrating \eqref{eq:derivu} from $t$ to $\infty$
shows that $\log u$ converges exponentially as well. Finally, passing to a limit in \eqref{eq:derivu}  shows
that the limiting metric has $R \equiv \rho$. 

\medskip

\noindent{$\bf{\rho = 0 :}$}
It is simplest in this case to proceed by writing $g(t)$ relative to a background flat metric, $g(t) = u(t) \bar{g}$ where $R(\bar{g}) \equiv 0$.
(Of course, in doing this we use uniformization.) Consequently, applying the maximum principle to the evolution equation 
$\del_t \phi = e^{-\phi}\Delta_{\bar{g}} \phi$ for $\phi = \log u$ one finds that $|\log u| \leq C$ for all $t \geq 0$. This bounds
the diameter, injectivity radius and Sobolev constants for $g(t)$ at all times.  It is then not hard, although somewhat
lengthy, to derive a sequence of bounds for all higher Sobolev norms of $u$ and show that all higher derivatives
tend to zero exponentially. Details appear in \cite{Ha1}. Hence $u \to \mbox{const.}$ and $g(t)$ converges to a flat metric. 

\medskip

\noindent{\textbf{$\rho > 0$ and $R(0) \geq 0$:}} 
To handle this case, Hamilton in \cite{Ha1} proves monotonicity of an ``entropy" function and introduces a generalization of the Li-Yau Harnack inequality \cite{LY}. 

To describe the Harnack inequality, consider a one-parameter family of Riemannian metrics $g(t)$ which is $\calC^1$ in $t$. 
For any two points $(x,T), (\xi,\tau) \in M^2 \times \RR^+$,  define
\begin{equation}
\label{Delta}
\Delta(\xi,\tau,x,T) = \inf_{\gamma}\int_{\tau}^T\left(\frac{ds}{dt}\right)^2\, dt,
\end{equation}
where the infimum is taken over all paths parametrized by $t \in (\tau,T)$ which connect $(\xi,\tau)$ to $(x,t)$;
here $\frac{ds}{dt}$ is the velocity in space at time $t$. 

\smallskip

{\bf Hamilton's Harnack estimate:} Let $g(t)$ be a solution to the Ricci flow equation on a compact surface $M^2$ with 
$R > 0$ for $0 < t \le T$. For any two points $(x_1,t_1)$ and $(x_2,t_2)$ in space-time with $0 < t_1 < t_2$ 
one has the inequality
\[
(e^{\rho t_1} - 1) R(x_1, t_1) \le e^{\Delta/4}(e^{\rho t_2} - 1) R(x_2,t_2),
\]
which controls the relative values of $R(x_1,t_1)$ and $R(x_2,t_2)$.

\smallskip

{\bf Entropy estimate:} For any  one-parameter family of  metrics $g(t)$ on $M^2$ with scalar curvature $R(t)$, define
\[
\mathcal{E}(t) := \int_{M^2} R\log R\, dA;
\]
this is the (negative of the) entropy of that family. If $g(t)$ is evolving via Ricci flow and if $R(t)>0$, then Hamilton
proves in \cite{Ha1} that $\mathcal{E}(t)$ is nonincreasing. 

\smallskip

The monotonicity of $\mathcal{E}(t)$ gives a uniform upper bound for $\int_{M^2} R\log R\, dA$; the Harnack estimate 
then can be used to show that $R$ itself is uniformly bounded above. A standard geometric estimate from \cite{CE} 
shows that the injectivity radius of $(M^2,g(t))$ is at least $\pi/\sqrt{R_{\max}(t)/2}$, hence 
is bounded away from zero. The same estimate also gives a uniform upper bound on the diameter. 

The diameter bound and Harnack estimate together show that for any $t \geq 1$ and any $x,y \in M^2$, we have
\[
R(x,t) \le CR(y,t+1),
\]
whence $0 < c_1 \le R \le c_2$ for all $t \ge 0$. 

Further control of the curvature is obtained as follows. Consider the evolution of the Hessian quantity, $|M_f|^2$:
\[
\del_t |M_f|^2 = \Delta_{g(t)} |M_f|^2 - 2|\nabla M_f|^2 - 2R|M_f|^2
\]
(here $f$ is the potential function for $R(t)$). Applying the maximum principle yields $|M_f| \le Ce^{-ct}$ for $c > 0$.

To complete the proof, Hamilton considers the modified Ricci flow equation
\[
\del_t g = (\rho - R)g - 2\nabla^2 f,
\]
the solutions of which differ from those of (\ref{eq:rfs0}) by diffeomorphisms. He shows that any solution of this modified
equation converges exponentially to a limiting metric with $M_f \equiv 0$. Metrics for which $M_f$ vanishes identically
are Ricci gradient solitons.\footnote{A Ricci soliton is a solution of the Ricci flow which evolves via diffeomorphism and
dilation; it is a Ricci gradient soliton if the vector field generating the diffeomorphisms is the gradient of a scalar function $f$.} 
The final step, therefore, is to establish the following:
\begin{thm}
\label{thm-sphere}
If $M^2$ is compact, the only gradient soliton metrics are those with constant curvature.
\end{thm}
This is due to Hamilton \cite{Ha1}, with an alternate proof using the Kazdan-Warner identity appearing in \cite{CK}.  
Both proofs appeal to uniformization.

\medskip

\noindent{\textbf{$\rho > 0$ and $R(0)$ changes sign:}}
\label{any sign}
The key new idea to handle this case is a modified entropy functional introduced by Chow \cite{Ch1}. 
Consider the PDE \eqref{eq-R1} satisfied by $R$, and the associated ODE $\frac{d}{dt}s(t) = s (s-\rho)$ 
obtained by dropping the Laplacian term. This ODE has the solution $ s(t) = \rho (1 - 
(1-\frac{\rho}{s_0})e^{\rho t})^{-1}$. The function $R-s$ satisfies
\[
\frac{\partial}{\partial t} (R - s) = \Delta(R-s) + (R-r+s)(R-s).
\]
Hence if one chooses $s_0 < R_{\min}(0) < 0$, then it follows from  the maximum principle that $R - s > 0$ on the entire interval of existence
of $g(t)$. 

Chow's {\it modified entropy} is the functional
\[
N(g(t),s(t)) = \int_{M^2} (R-s)\log(R-s)\, dA.
\]
Although this is not monotonic, it does remain uniformly bounded. There is also a modified Harnack inequality:
\[
R(x_2,t_2) - s(t_2) \geq e^{-\Delta/4 - C(t_2 - t_1)}(R(x_1,t_1) - s(t_1)),
\]
again for any $(x_1,t_1), (x_2,t_2) \in M^2 \times \RR^+$ with $t_1 < t_2$, and with $\Delta$ as in \eqref{Delta}. 

This Harnack estimate and the upper bound on $N(g(t),s(t))$ yield a uniform upper bound on curvature, and from
this we obtain a uniform lower bound on the injectivity radius and a uniform upper bound on the diameter. 
The diameter bound and Harnack estimate now give $R - s \ge c > 0$ for all $t \geq 0$, but since $s(t) \to 0$ 
exponentially, we see that $R(t)$ eventually becomes strictly positive when $t$ is large enough. This reduces
us to the previous case, and completes the proof of convergence to a constant curvature metric in all cases. 

\subsection{Uniformization theorem}
As  pointed out explicitly in the above discussion of the proofs of Theorems \ref{thm-ham} and \ref{thm-chow}, we have  invoked the uniformization theorem in the case $\rho = 0$. It turns out that it is not 
too difficult to circumvent this and carry out the argument for that case using an 
arbitrary background metric. However, uniformization is also used in a more subtle way when $\rho > 0$ in the proof of 
the monotonicity of the entropy and in the proof of Theorem \ref{thm-sphere}. Hence the proof sketched above does not give a 
new proof of uniformization for surfaces. However, with more work, one is able to establish the uniformization
theorem using Ricci flow after all.  In particular, a different proof of the entropy estimate which does not presuppose uniformization 
appears in \cite{Ch2}; and \cite{CLT} contains an alternate proof of Theorem \ref{thm-sphere}
which also does not invoke uniformization, but involves the following lemma as the  key ingredient:
\begin{lem}
\label{lem-rot-symm}
Let $(M^2,g)$ be a complete Riemannian surface with a nontrivial Killing vector field $X$. If $X$ vanishes at a point 
$q\in M^2$, then $(M^2,g)$ is rotationally symmetric about that point. 
\end{lem}

To deduce Theorem \ref{thm-sphere} from this lemma, one argues as follows. If $g$ is a shrinking Ricci gradient soliton, then the equations
\[
\Ric = g + \nabla^2 f \Longleftrightarrow (R/2 - 1)g = \nabla^2 f
\]
imply that $\nabla f$ is a conformal Killing field. If $J$ is an almost complex structure on $M^2$, then it follows from a 
short computation that $J(\nabla f)$ is actually Killing. Since $M^2$ is compact and two dimensional, $\nabla f$ (and hence 
$J(\nabla f)$) must vanish somewhere, so by Lemma \ref{lem-rot-symm}, $(M^2,g)$ is rotationally symmetric and one can write
$g = dr^2 + h(r)^2\, d\theta^2$, for some nonnegative function $h$. If we now substitute this form of the metric into 
the soliton equation, we obtain a system of ODEs in $h$ and $f$. It is not hard to solve this system, and by
direct inspection of the solution, one sees that $g$ has constant curvature. 

\section{Open surfaces}
The behaviour of Ricci flow for open surfaces is less well understood than for closed surfaces. Indeed, besides the questions of 
long time existence and convergence, a number of new issues arise for open surfaces, including the nontriviality of proving short-time 
existence and the possible nonuniqueness of the solution for a given initial metric. It is also of interest to know whether $g(t)$ 
remains quasi-isometric to $g_0$ (i.e., if $c_1 g_0  \leq g(t) \leq c_2 g_0$ for some positive constants $c_1 < c_2$).

There has been recent significant progress toward understanding these questions.  It has been known for
some time \cite{Shi} that if $(M^2,g_0)$ is complete and has curvature bounded above, then there is a unique solution with initial 
data $g_0$ on some small time interval. If on the other hand $(M^2,g_0)$ is incomplete, this may no longer be true, 
and we discuss some instances of this below.  Amongst all possible solutions with a given $g_0$ there is one which 
is `maximally stretched' and has the remarkable property that it becomes instantaneously complete, i.e., $g(t)$ is complete 
for any $t > 0$, regardless of the geometric nature of $g_0$. This was discovered by Topping \cite{T}, and later elaborated 
upon in much greater detail by him and Giesen \cite{GT1}, \cite{GT2}. 

As to whether $g(t)$ remains quasi-isometric to $g_0$, we restrict attention to surfaces with finite topology and recall 
the classical fact that the nature of a uniformizing metric in any conformal class of such an open surface depends both 
on the Euler characteristic and on the conformal type of the ends, i.e., whether they are conformal to annuli or to punctured disks. 
(These are the only two possibilities.)  If an end is conformal to an annulus, then the uniformizing metric is hyperbolic with 
infinite area, and the metric on that end is a `hyperbolic funnel',  of the form $dr^2 + e^{2r} d\theta^2$, $\theta \in S^1_\ell = 
\RR/ \ell \ZZ$. On an end which is conformal to a punctured disk, if the metric is hyperbolic, then it is a hyperbolic cusp, 
taking the form $dr^2 + e^{-2r}d\theta^2$. Thus the only (orientable) complete uniformizing metrics are the following:  the 
plane $\RR^2$; the hyperbolic plane $\HH^2$; the hyperbolic cylinder (which is a funnel on one end and a cusp on the other); the flat 
cylinder $\RR \times S^1$;  and then any complete hyperbolic surface (this is true even without assuming finite topology).  
If the asymptotic geometry 
of $(M^2,g_0)$ is different from that of its conformally uniformizing metric, then it is of interest to determine whether the change in
quasi-isometry type occurs at $t=0$, or at some later finite time initiated by reaching some critical stage of the geometry
in the flow, or only in the limit as $t \nearrow \infty$. This last case is the most plausible, and is now known to occur
for surfaces which are asymptotically conic \cite{IMS}, for example. However  the other possible behaviors have not been ruled out yet for other classes of complete initial metrics. 

One further issue which complicates the study of the limiting behaviour of the flow is the existence of nontrivial complete Ricci solitons.
The main example is the so-called cigar soliton, which is given by the explicit expression (in terms of the standard Euclidean coordinates $(x,y)$ on $\RR^2$)
\begin{equation}
\label{cigar}
g_{\mathrm{sol}}= \frac{dx^2 + dy^2}{1 + x^2 + y^2}. 
\end{equation}
It is known that this is the only soliton which exists for all $t \in \RR$ and such that each $g(t)$ has bounded
curvature; there are other solitons which exist on half-lines $(-\infty, T)$ or $(T,\infty)$, but a complete classification
of complete solitons in two dimensions has not yet been attained. 

Based on all of this, the natural problem for complete surfaces is to determine whether the flow
converges (possibly jumping quasi-isometry class) to the complete conformal uniformizing metric
or alternately to a soliton. We survey the state of knowledge on this. For incomplete metrics, on the other
hand, the natural problem is to understand ways that one might impose `boundary conditions' at the points 
of $\overline{M^2}\setminus M^2$, and how these can determine different choices of the flow.  

\subsection{Different classes of complete surfaces}
We now discuss the detailed results that have been obtained for various classes of open surfaces with 
specified conditions on their asymptotic geometries. This is divided into the study of surfaces
conformal to $\RR^2$, and those for which the uniformizing metric is hyperbolic. 

\subsubsection{Surfaces conformal to $\RR^2$}
Suppose that $M^2=\RR^2$, and that the metric $g_0$ is conformal to the standard Euclidean metric, so using the complex 
notation $z=x_1+ix_2$ and $|dz|^2 = (dx_1)^2+(dx_2)^2$, we have $g_0 = u_0 |dz|^2$. The expectation is that, under the 
flow \eqref{eq:rfs0} with $\rho = 0$, at least some classes of metrics which are conformal to the Euclidean metric, or even just 
quasi-isometric to it,  converge to $|dz|^2$ as $t \to \infty$. This is indeed true, as we discuss below. However, other behavior 
can also occur; for example,  the cigar soliton \eqref{cigar} is conformally Euclidean and does not flow to a flat metric.
Thus it is necessary to find some way to separate out those metrics whose flow might converge to either the flat metric
or cigar soliton, or perhaps even to some alternative.  This is accomplished using the following two notions:

\begin{defn}
Let $(\RR^2, u |dz|^2)$ be a complete metric on the plane. The \emph{circumference at infinity} of $g$ is defined as
\begin{multline*}
C_\infty(g)= \\
\sup_K \inf_{D} \{L(\partial D)\ |\ \forall~ \mbox{compact sets}~K\subset \RR^2 \ \mbox{and all open sets}~ D \supset K \};
\end{multline*}
similarly, the \emph{aperture} of $g$ is defined by
\[
A(g)=\lim_{r \rightarrow \infty} {{L(\partial B_r)} \over {2\pi r}},
\]
where $B_r$ denotes the geodesic ball of radius $r$ and $L(\partial B_r)$ is the length of its boundary.
\end{defn}

These measure slightly different things, but clearly $(\RR^2, |dz|^2)$ has infinite circumference at infinity
and aperture equal to $1$, while  $(\RR^2, g_{\mathrm{sol}})$ has circumference at infinity equal to $2\pi$ 
and aperture $0$. On the other hand, a metric which pinches off to a cusp has both its circumference at
infinity and its aperture equal to $0$.  Wu shows that if $(M,g(t))$ has finite total curvature for each $t$,
then the aperture $A(g(t))$ is independent of $t$. 

\medskip

\noindent \textbf{Conditions for convergence to the Euclidean metric}
\label{WuIJ}

The first results to address these questions on open surfaces were obtained by Wu \cite{W2}, soon  after the
initial work of Hamilton for the closed case.  To state her result, we introduce a slightly weaker notion of convergence.
A solution is said to have \emph{modified subsequence convergence} if there exists a one-parameter family of 
diffeomorphisms $\{\varphi_t\}$ such that for any sequence $t_j \rightarrow \infty$, there exists a subsequence 
(denoted again by $t_j$) for which $\varphi_{t_j}^*g_{t_j}$ converges uniformly on every compact set as $j \rightarrow \infty$. 

\begin{thm} \label{wu} \cite{W2}
Let $g(t)=u(t,z)|dz|^2$ be a solution to \eqref{eq:rfs0} with $\rho=0$ such that $g(0)$ is a complete metric with bounded curvature 
and $|\nabla \log u_0|$ is uniformly bounded. Then $g(t)$ converges in the modified subsequential sense as $t\rightarrow \infty$. 
Moreover if $g(0)$ has positive curvature, then any metric which arises in such a limit is necessarily either the cigar soliton \eqref{cigar} 
if $C_\infty(g(0))<\infty$, or else the flat metric $|dz|^2$ if $A(g(0))>0$.  
\end{thm}

The proof relies mainly on the maximum principle.  Wu uses  Shi's general results \cite{Shi}, which hold in all dimensions, 
to obtain short time existence and uniqueness; this uses the assumption on bounded initial curvature.  The strategy for long-time 
existence is somewhat like that  in the compact case, but involves the quantity $h:= R+ |\nabla \log u |^2$, which satisfies the 
evolution inequality $\del_t h \le \Delta h$.  The maximum principle then implies bounds on $R$ and $|\nabla \log u|$ on any
finite time interval $[0,T)$.  Applying similar arguments to analogues of $h$ which involve higher derivatives of $\log u$ and $R$, 
one obtains uniform bounds on all of these as well. These a priori bounds yield long-time existence.

To show that $u$ converges, note that \eqref{eq-R1}  implies that if $R(0) > 0$, then $R(t) > 0$ for all $t > 0$. The equation $\del_t u = -Ru$
implies that $u$ decreases monotonically. Since $u(t) > 0$ for all $t$, it must converge pointwise; the bounds on $|\nabla \log u|$ then
give uniform convergence on any compact set. 

Convergence of the conformal factor $u$ does not guarantee convergence of $g(t)$ to a Riemannian metric. Indeed, $u$ is positive,
but could converge to zero in some regions or even everywhere; this happens, for example, in the flow corresponding to the cigar 
soliton which has evolving conformal factor $u(x,y.t)= (e^{4t} +x^2 + y^2)^{-1}$. This can be `fixed' using a family of 
diffeomorphisms, to obtain modified subsequential convergence. The proof of Wu's theorem is completed by examining the geometric 
implications of the two conditions $C_{\infty}(g) < \infty$ and $A(g) > 0$.

One would like  to remove certain hypotheses in this theorem, such as  the bounds on curvature and on $|\nabla \log u_0|$. This 
is partially accomplished in the recent paper \cite{IJ} by the first author and Javaheri, in which the condition that the initial 
curvature be positive is replaced by the assumption that $u_0$ is bounded.

\begin{thm} \label{IJ} \cite{IJ}
Suppose that $g_0 = u_0 |dz|^2$ has bounded curvature and that $\log  u_0 \in \calC^1(\RR^2) \cap L^\infty$. Then the solution
$g(t)$ exists for all $t \ge 0$, and for each $k \in {\mathbb N}$ it has modified subsequential convergence in $\calC^k$ to $|dz|^2$.
\end{thm}

The bound on  $\log u_0$ shows that $g_0$ is complete and has infinite area (hence it avoids singularity formation at finite time,
as discussed below). The bounds on curvature and standard elliptic estimates yield that $|\nabla \log u_0| \leq C$. 
Much as in the proof of Theorem \ref{wu}, a considerable amount  of the analysis here involves a succession of maximum 
principle arguments controlling the higher derivatives. Note that the boundedness of $\log u_0$ implies that $A(g_0) > 0$, 
so once uniform convergence is proved, the limit must be $|dz|^2$ rather than $g_{\mathrm{sol}}$.  If the condition 
$|\log u_0| \leq C$ is replaced by certain bounds in terms of the cigar soliton's conformal factor $1/(1+x^2 + y^2)$, then it is 
possible to show that the corresponding solution $g(t)$ converges to $g_{\mathrm{sol}}$. 

The recent paper \cite{GT2} contains a slight improvement for  this case.
\medskip

\noindent \textbf{Finite time singularities}
\label{R2Sing}

As discussed in \S \ref{GenCon},  for metrics on closed surfaces, the average of $R$ over $M^2$ is a good choice for $\rho$ since 
the corresponding
flow fixes the area of $(M^2,g(t))$.  For certain complete surfaces one might hope for a similar effect, and as we discuss
below, this is more or less correct for surfaces with finite total curvature. Conversely, if $(M^2,g_0)$ has finite
area, then the unnormalized flow with $\rho = 0$ does not preserve the area, and a similar calculation to the one at
the end of \S 2 should lead to the formation of a singularity at time $T_0 = (1/4\pi)\mathrm{Area}\,(M^2,g_0)$. 
However, the flow is in fact more singular than this, since it turns out that not only does the flow not necessarily have unique
solutions, but many of these solutions develop singularities earlier than the `expected' time above.  This phenomenon
was discovered by Daskalopoulos and del Pino \cite{DD1}. We describe their findings in a bit more detail now.

Suppose that $g_0 = u_0 |dz|^2$ and in addition that
\[
\mbox{Area}\,(M,g_0) = \int_{\RR^2} u_0\, dx < \infty.
\]
In fact, these authors are not necessarily thinking of the geometric problem, but rather considering  general solutions of
\eqref{eq:rfs} with $R_0 = \rho = 0$ and $g_0 = |dz|^2$. Hence they also allow initial data $u_0$ which are
nonnegative, possibly even with compact support, and with finite integral. They prove that for any 
$\gamma  \geq 2 $, there exists a  solution $u_{\gamma}(x,t)$ of \eqref{eq:rfs} (with $R_0 = \rho = 0$) such that
\begin{equation}\label{eqn-intc}
\int_{\R^2} u_\gamma (x,t) \, dx = \int_{\R^2} u_0 \, dx  -  2\pi \gamma \, t,
\end{equation}
and then observe that this solution exists only on the interval $0 \leq t < T_\gamma := 
(1/2\pi\, \gamma) \, \int_{\R^2} u_0 \, dx.$ 

They also consider maximal solutions of this same problem. The solution $u$ is called \emph{maximal} if, for any other solution 
$\tilde{u}$ with the same initial data as $u$, one has $\tilde{u}(x,t) \leq u(x,t)$ for all $t \geq 0$. The maximal solutions 
are simply the solutions $u_2$ described above, and these develop singularities at time $T_2 = (1/4\pi)\int u_0  dx$.
Moreover, the maximal solution $u_2$  for a given set of initial data is complete and is unique in the category of maximal 
(complete) solutions. It is shown in \cite{DH} that these maximal solutions develop type $II$ singularities at this 
extinction time $T_2$ in the sense that 
\begin{equation}\label{eqn-c}
\frac{c}{(T_2-t)^2} \leq R_{\max}(t) \leq \frac{C}{(T_2-t)^2}  
\end{equation} 
for constants $0 < c < C$ and all $0< t < T_2$. The papers \cite{DD2} (in the rotationally symmetric case) 
and \cite{DS} (in general) rigorously verify the formal analysis of King \cite{K} and give precise asymptotics 
for this type $II$ singularity formation. The final result is as follows:

\begin{thm}(Daskalopoulos, del Pino, Sesum)
Let $u_0 \geq 0$ have finite integral and let $u(x,t)$ be a maximal solution to \eqref{eq:rfs} with initial data $u_0$. 
Then $u$ becomes extinct at time $T_2 < \infty$. If $u_0$ is compactly supported,  then the flow behaves differently
in the two regions $(T_2-t)\, \log  r > T_2$ and $(T_2-t)\, \log r \leq  T_2$. In  the first, called the {\em outer region}, 
after an appropriate change of variables and rescaling, the solution converges to a hyperbolic cusp solution, 
$2 t\, /|x|^2 \, \log^2 |x|$ (corresponding to the flow of the hyperbolic cusp on $\RR^2\setminus \{0\}$).
In the second {\em inner region}, the solution decays exponentially and converges, after an appropriate change of 
variables and rescaling, to a soliton, which in this context is the solution $U(x,\tau) = 1/( \lambda  |x|^2 + e^{4\lambda \tau})$, 
with  $\tau=1/(T_2-t)$ and $\lambda=T_2/2$, of the equation $U_\tau = \Delta \log U$ for $\tau \in \RR$. 
\end{thm}

\subsubsection{Surfaces conformally covered by $\HH^2$}
If the complete open surface $(M^2,g_0)$ is not conformally covered by $\RR^2$, then we can identify $M^2$ as a 
quotient $\HH^2/\Gamma$, where $\Gamma$ is a discrete subgroup of the group $SO(2,1)$ of conformal transformations of the disk, and where the metric $g_0$ lifts to a $\Gamma$-invariant metric $\tilde{g}_0$ on the unit 
disk $\mathbb{D}$ which is conformal to the hyperbolic metric $g_{\mathrm{hyp}}$: $\tilde{g}_0 = u_0 \, g_{\mathrm{hyp}}$. 
In this case the Ricci evolution $g(t)$ of $g_0$ on $M^2$ lifts to the solution $u(z,t)\, g_{\mathrm{hyp}}$ on $\mathbb{D}$,
where $u$ is a solution of \eqref{eq:rfs} (relative to $g_{\mathrm{hyp}}$) with initial condition $u_0$. In particular,
a complete solution of \eqref{eq:rfs0} on $\mathbb{D}$ with complete $\Gamma$-invariant initial metric 
$\tilde{g}_0$ remains $\Gamma$-invariant; see \cite{GT2}. In principle, the study of the Ricci evolutions of
$g_0$ on $M^2$ or $\tilde{g}_0$ on $\mathbb{D}$ are equivalent. In practice, however, the former is more useful if one 
is interested in keeping track of the asymptotics of each of the metrics $g(t)$ on $M$, while the latter method is 
simpler in certain senses since one does not need to deal directly with the topological or geometric complexity of $M^2$. 

\medskip

\noindent \textbf{Complete metrics on the disk}

The first work in this setting was done by Angenent, who contributed an appendix to Wu's paper \cite{W2}, in which he presents 
a few basic results. His primary observation is that the solution to \eqref{eq:rfs0} (with $\rho=0$) with 
initial condition $g_{\mathrm{hyp}}$ is simply $g(t) = (1+2t)g_{\mathrm{hyp}}$, and that constant multiples of this can be used as either 
an upper or lower barrier for the evolution of metrics which are quasi-isometric to the hyperbolic metric.

\begin{prop}[Angenent \cite{W2}]
Let $g_0$ be a complete metric on $\mathbb{D}$ which satisfies $c_1 g_{\mathrm{hyp}} \leq g_0 \leq c_2 g_{\mathrm{hyp}}$ 
for some constants $0 < c_1 < c_2$. Then $c_1 (1+2t) g_{\mathrm{hyp}} \leq g(t) \leq c_2(1+2t) g_{\mathrm{hyp}}$,
and hence in particular, $g(t)$ exists for all time.
\end{prop}
Note that this hypothesis holds if $g_0$ is the lift of a metric on a compact quotient of the disk.  On the other hand, 
Angenent also notes that the evolution of $(1-|z|^2)^{-k}$ can also be used as a barrier for the same equation for any $k \in \mathbb{N}$, 
so that there are many other complete families of solutions which are not equivalent to the hyperbolic metric
at any time. 

\medskip

\noindent \textbf{Surfaces with asymptotically hyperbolic cusp ends} 

The alternate point of view, to study the Ricci flow directly on complete open surfaces, has only been undertaken
more recently.  The first case to be studied was that of finite area surfaces with ends which are asymptotic
to hyperbolic cusps. This was done by Ji together with  the second and third authors of this paper \cite{JMS}.

To describe this result, suppose that $(M^2, g_0)$ is a complete surface with finite area such that $M^2$ has
only finitely many ends, and each end $E_j$ is diffeomorphically identified with a half-cylinder 
$\RR^+ \times S^1$, with coordinates $(s,y)$. We also assume that for each end, 
\[
\left. g_0 \right|_{E_j} = a_j(ds^2 + e^{-2s}dy^2) + k,
\]
where $|k| = \calO(s^{-\tau})$ for some $\tau > 1$. The pointwise norm of the symmetric $2$-tensor $k$ is
taken with respect to the model hyperbolic cusp metric $ds^2 + e^{-2s}dy^2$, and the actual hypothesis
is that $k$ lies in some weighted H\"older space with weight factor $s^{-\tau}$; the coefficients $a_j$
are positive constants, and their inclusion allows us to assume that $g_0$ has different limiting 
curvatures along each end. 

There are, as usual, three issues to consider: short-time existence, long-time existence and convergence. 
Perhaps surprisingly, it is the first two of these that require more effort here, but once these are established,
the convergence of the solution $g(t)$ to a complete finite area hyperbolic metric proceeds almost exactly
as in the compact case.

If one is interested in understanding the spatial asymptotics of each metric $g(t)$, then it is not sufficient to 
quote \cite{Shi} for short-time existence of the Ricci flow on complete manifolds of bounded geometry. 
For simply proving that the solution exists, 
this may not seem important, but if one is interested in understanding the evolution of any nonlocal
geometric or analytic quantities, (for example, the renormalized determinants considered in \cite{AAR}), then
the fact that $g(t)$ has asymptotics of a fixed type for each $t$ is very important.  To handle this, it is shown in
\cite{JMS} that the linearization of \eqref{eq:rfs} preserves the weighted H\"older spaces which contain the initial 
metric $g_0$, and hence the standard contraction mapping argument then applies to guarantee short-time
existence in this class of metrics.  One further interesting feature is that one must understand how this
short-time flow alters the asymptotic curvatures (i.e., the constants $a_j$) on each end.

In this setting, it is possible to choose the constant $\rho$ in \eqref{eq:rfs0} to preserve the area. The correct
choice is to set $\rho$ equal to the average of the scalar curvature $R$, but one must check that this average
is independent of $t$; this, in turn, is equivalent to the fact that metrics of this type satisfy a Gauss-Bonnet theorem.
The proof of long-time existence closely follows the argument in the compact case, but the difficulty is in showing
that there exists a potential function $f$; i.e., a solution to $\Delta f = R - \rho$, which has bounded gradient. 
This equation is solvable by standard $L^2$ methods since the right hand side has vanishing integral and $0$ is isolated
in the spectrum of $\Delta$, but it is not clear that this solution has bounded derivative. 
The alternate argument, to show that there is a solution with
these properties, is the heart of that paper, and it involves a detour into the scattering theory
of such metrics inasmuch as one needs to determine the possible asymptotics at infinity of solutions to
equations of the form $\Delta h = q$ for $q \in \calC^\infty_0(M)$. 

Having obtained a suitable potential function, one then proceeds exactly as in the compact case to show
that $R(z,t)$ converges to a negative constant and that the conformal factor $u(z,t)$ converges. 

The final result is the following:
\begin{thm} \cite{JMS}
\label{thm-JMS}
Let $(M^2,g_0)$ be a complete surface with finite area, with $\chi(M^2) < 0$, and with  $g_0$ satisfying  the asymptotic 
conditions described above. Then the flow (\ref{eq:rfs0}) with $\rho$ equal to the average of the scalar curvature
of $g_0$  has a unique solution which exists for all $t > 0$ and which converges exponentially to some multiple
of the unique complete hyperbolic metric in the conformal class of $g_0$. 
\end{thm}

\noindent \textbf{Surfaces with infinite area asymptotically hyperbolic ends}

Another case very closely related to the one just discussed is that in which  $(M^2,g_0)$ has ends which
are asymptotic to either hyperbolic cusps or hyperbolic funnels (i.e., to metrics of the form
$ds^2 + e^{\pm 2s} dy^2$). The reason that this is conceptually similar is that the expected 
limiting metric is the uniformizing hyperbolic metric, which as per the discussion at the beginning
of this section has the same quasi-isometry type as $g_0$. 

The proof of existence and convergence in this setting has been handled in another recent paper by 
Albin, Aldana and Rochon \cite{AAR}. Their aim in this paper is not only to study the convergence
of the Ricci flow, but also to study the behaviour of the renormalized determinant under this evolution.
They generalize the well-known result by Osgood-Phillips-Sarnak \cite{OPS}, who study
the gradient flow of the determinant functional in conformal classes on compact Riemann surfaces,
and show that the determinant is monotonic under this Ricci flow (This had been shown on compact
surfaces independently by M\"uller and Wendland \cite{MW} and by Kokotov and Korotkin \cite{KK}). 

In any case, their main theorem about the Ricci flow itself is the following:
\begin{thm} \cite{AAR}
\label{thm-AAR}
Let $(M^2,g_0)$ be a complete surface with infinite area and asymptotically hyperbolic cusp or funnel ends. 
Then the flow (\ref{eq:rfs0}) with any constant $\rho < 0$ exists for all $t > 0$ and converges exponentially 
to the unique complete hyperbolic metric in the conformal class of $g_0$. 
\end{thm}

The same ingredients are needed: short-time existence in the class of metrics with this specified
asymptotic type, long term existence via Hamilton's method involving potential functions and convergence.
In this case, when none of the ends are cusps, the potential function is much easier to obtain. Note
too that one does not need to make any assumption about $\chi(M^2)$.  

\medskip
\medskip

\noindent \textbf{Surfaces with asymptotically conical ends}

The final type of asymptotic geometry which has been studied is that of asymptotically conical (or AC) surfaces
$(M^2,g_0)$. The definition of such a surface is similar to the cases above: $M^2$ is required to have finite
topology, with each end $E_j$ identified with a half-cylinder $\RR^+ \times S^1$, and with the metric 
\begin{equation}
\left. g \right|_{E_j} = dr^2 + \alpha^2 r^2 dy^2 + k, \qquad \mbox{where}\ |k| = \calO(r^{-\tau}).
\label{acmetric}
\end{equation}
The constant $\alpha$ determines the cone angle $2\pi \alpha$ at infinity, and is the same as the aperture
we defined earlier. As before, the norm of the tensor $k$ is taken with respect to the model conic metric,
and the actual hypothesis is that $k$ lies in a certain weighted H\"older space.

An AC surface has finite total curvature; hence by a classic theorem of Huber is conformally finite, so each
end is conformal to a punctured disk. Provided $\chi(M^2) < 0$, the uniformizing metric is a hyperbolic metric 
with finite area and cusp ends, which is quite different than the initial metric $g_0$. It is interesting 
to determine how this transition happens. This is the topic of a recent paper by the three authors of
the present paper: 

\begin{thm}[\cite{IMS}] 
\label{thm-IMS}
Let $(M^2,g_0)$ be a surface with asymptotically conical ends and $\chi(M^2) < 0$. Then the Ricci flow with initial
metric $g_0$ exists for all time $t > 0$, and for each $t$, $g(t)$ is again asymptotically conical. There exists
an absolute constant $C_2$, and for any compact set $K \subset M^2$ another constant $C_1(K)$, such that
\[
C_1(K) (1+t) g_0 \leq g(t) \leq C_2 (1+t) g_0.
\]
The rescaled metric $\tilde{g}(t) := t^{-1} g(t)$ is thus locally uniformly bounded. It converges to the unique 
complete uniformizing hyperbolic metric in the conformal class of $g_0$. 
\end{thm} 

Let us make a few comments about the proof, which is rather different from the proofs of  the other two cases discussed above.
One starts by proving short-time existence in the class of AC metrics; this is straightforward. The existence
of the potential function with bounded gradient is proved using some well-known facts about the mapping
property of the Laplacian on AC spaces; this is a simpler version of the corresponding fact in the asymptotically
hyperbolic cusp case.  These two steps give long-time existence.  The final step, however, has a different flavor than
the arguments sketched above (and in the compact case).  The first key point is that the scaled metric $\tilde{g}(t)$ 
satisfies a Ricci flow type equation, $\partial_t \tilde{g}(t) = -(1 + \tilde{R}) \tilde{g}$, where $\tilde{R}$ is the scalar
curvature of $\tilde{g}$. A maximum principle shows that the conformal factor $\tilde{u}(t,z) = u(t,z)/t$ is monotonically
decreasing, hence converges pointwise. A succession of compactness arguments (along with the condition $\chi(M^2) < 0$)  prove that the limiting
function is smooth and vanishes nowhere, and from this we conclude easily that the limiting metric $\tilde{g}_\infty$
is hyperbolic. 

\section{Flows on incomplete surfaces}

The previous discussion has focussed on the Ricci flow for complete metrics. If $(M^2,g_0)$ is not complete, 
then the flow may not be unique. We discuss two aspects of this problem, first considering solutions which
are complete for all $t > 0$ and then considering solutions which remain incomplete.

\subsection{The instantaneously complete flow}
It is natural to search for solutions which are maximal in some sense, since for these, one might hope to prove uniqueness. 
An important new direction has been opened with the recent notion of an instantaneously complete solution, 
appearing originally in the work of DiBenedetto and Diller \cite{DiDi}, but taken up in a more general way by Topping 
\cite{T} and then simplified and generalized in two papers by Giesen and Topping \cite{GT1}, \cite{GT2}. 
Their results are both striking and strong, and subsume much of what has been discussed above. We describe 
some of their ideas and results now.

Suppose that $(M^2,g_0)$ is any open surface with a not necessarily complete metric. A solution $(M^2,g(t))$ of \eqref{eq:rfs0} 
is called {\em instantaneously complete} if it is a solution to this equation on some time interval $(0,T)$, if $g(t)$ is complete for 
every $t$ in this open interval, and if $g(t) \to g_0$ in some reasonable geometric sense (e.g., locally smoothly
on the regular set of $g_0$) as $t \searrow 0$. This solution is called {\em maximally stretched} if any other solution $\tilde{g}(t)$
of this flow with $\tilde{g}(0) = g_0$ satisfies $\tilde{g}(t) \leq g(t)$ on the common interval of existence.

The culminating and most recent result by Giesen and Topping is the following:
\begin{thm}[\cite{GT2}]
\label{GT-thm}
Let $(M^2,g_0)$ be any open surface. Then there exists a unique maximally stretched instantaneously complete
solution of \eqref{eq:rfs0} with initial condition $g_0$ which is defined on some interval $[0,T)$. The final
time $T$ is infinite unless $M$ is conformally covered by $\RR^2$ and the area of $(M^2,g_0)$ is finite, and in
this latter case, the area of $(M^2,g(t))$ tends to $0$ as $t \nearrow T$. If $(M^2,g_0)$ is conformally covered
by the hyperbolic disk, then $(2t)^{-1}g(t)$ converges locally uniformly to the hyperbolic uniformizing metric on $M$.
If $g_0 \leq C g_{\mathrm{hyp}}$, then this convergence is uniform on all of $M^2$ and decays in any
$\calC^k$ topology like $t^{-\gamma}$ for some $0 < \gamma < 1$. This solution is unique amongst
solutions with all of these properties. 
\end{thm}

They also consider the question of whether the instantaneously complete flow is unique if one drops the hypothesis
that $g(t)$ be maximally stretched and prove that this is indeed the case if $(M^2,g_0)$ is conformally covered
by the plane, and also if $(M^2,g_0)$ is conformally covered by the hyperbolic disk provided $g_0$ is bounded
above by some multiple of the hyperbolic metric. 

This set of results is a tour de force. Since its hypotheses allow $g_0$ to be complete, Theorem \ref{GT-thm} contains
each of the previous ones for the different classes of complete surfaces. The proof of this theorem relies
on a clever combination of comparison arguments, including Yau's generalization \cite{Y} of the classical Schwarz-Pick-Ahlfors
theorem for conformal self-maps of the disk. The main difficulty in proving uniqueness is that there is no direct assumption
on the behaviour of any competing Ricci  flow geometries near spatial infinity. This control is built up by 
appealing to geometric results which exploit completeness and a direct analysis of the conformal factor of the flows, more 
in the spirit of  the logarithmic fast diffusion equation. The second paper \cite{GT2} is an improvement over \cite{GT1}
in that the authors drop the assumption on the upper curvature bound of the initial metric, they give a precise formula for the
maximal existence time in all cases, and they show that the (rescaled) Ricci flow converges to a hyperbolic metric whenever 
$M$ is conformally covered by the disk. 

\subsection{Solutions which remain incomplete}
This final section contains a brief discussion of two further results which indicate the existence of multiple
solutions of \eqref{eq:rfs0} starting from an incomplete surface $(M^2,g_0)$ which all remain incomplete, and which 
correspond to different choices of `boundary conditions' at the frontier $\overline{M^2} \setminus M^2$.
The first case is for surfaces with boundary, with  the use of the term `boundary condition' being the
standard one for a parabolic problem; the second case is for surfaces with isolated conic singularities,
and in this case, by boundary conditions we mean particular choices of self-adjoint extensions of 
the Laplacian.  It is not clear how one might formulate boundary conditions for more general
incomplete metrics. 

\medskip

\noindent \textbf{Surfaces with boundary}

The analysis of Ricci flow for surfaces with boundary, for two different geometrically natural
boundary conditions, was carried out by Brendle a decade ago \cite{Br}. He first studies the normalized
Ricci flow with the condition that the boundary be totally geodesic. This is equivalent to imposing
a Neumann condition on the conformal factor $u$ in \eqref{eq:rfs}.   Brendle shows that, as expected,
this flow converges exponentially to a metric with constant Gauss curvature and geodesic boundary.
The second flow he studies requires that $g(t)$ remain flat for all $t > 0$; the boundary condition which
drives this flow is a dynamic one:
\[
\left. \del_t g \right|_{\del M^2} = \left. -2 (\kappa - \overline{\kappa}) g \right|_{\del M^2}, 
\]
where $\kappa$ is the geodesic curvature of the boundary and $\overline{\kappa}$ is its average. In this
case, the solution converges exponentially to a flat metric with boundary of constant geodesic
curvature. These results reprove and generalize the work of Osgood-Phillips-Sarnak \cite{OPS}, who proved
the existence of these metrics in some cases by elliptic methods. 

Brendle's proof proceeds by a series of a priori integral estimates; it would be interesting to know whether these
results could also be proved using more direct maximum principle arguments.

\medskip

\noindent \textbf{Surfaces with conic points}

The form of a metric on a surface near an isolated conic point is given in polar coordinates
$(r,y) \in [0,r_0) \times S^1$ by
\[
g_0 = dr^2 + \alpha^2 r^2 dy^2 + k,
\]
where $\alpha$ is a positive constant which determines the cone angle $2\pi \alpha$, and where the
tensor $k$ decays (relative to the model metric $dr^2 + \alpha^2 r^2 dy^2$) in a H\"older sense
like $r^\tau$ for some $\tau > 1$.  Note that this is very similar to the form of an AC metric near
the large end of a cone, but of course here the conic point $r=0$ is at finite distance. 
Surfaces with all cone angles less than $2\pi$ are special in several ways; for example,
if $(M^2, g_0)$ is a surface with isolated conic singularities, each with cone angle less than $2\pi$,
then for any two points $q_1, q_2 \in M$ (these may be conic points), there exists a length minimizing
geodesic connecting these two points which contains no conic points in its interior. This angle
condition is also crucial for verifying certain properties of the Ricci flow on such surfaces.

One can also express a conic metric in conformal form,  
\[
g_0 = e^{2\phi_0} |z|^{2(\alpha-1)} |dz|^2 = e^{2\phi_0 + 2(\alpha-1)\log |z|} |dz|^2.
\]
Here $\phi_0$ is a bounded function which near $r=0$ has the form $\phi_0 = a + \tilde{\phi}_0$, where $\tilde{\phi}_0$
decays (in a H\"older sense) like $r^\tau$ for some $\tau > 0$. This makes it clear that the conformal structure 
determined by $g_0$ extends across the conic point. It also shows that altering the conformal factor by a multiple
of $\log |z|$ changes the cone angle. We are making a distinction between $|z|$ and the polar distance $r$, but 
a quick calculation shows that $r \approx c |z|^{\alpha}$, so $\log r \approx c' \log |z|$. 

One immediate difficulty is that a general metric of this form has unbounded curvature. Suppose, for example, 
that we write $g_0 = dr^2 + r^2h(r) dy^2$ (this is a third equivalent way to write any conic metric), where we suppose 
for simplicity momentarily that $h \in \calC^\infty([0,1))$ with $h(0) = \alpha^2$. Then the limit of the Gauss curvature 
of $g_0$ at $p$ blows up like a multiple of $h'(0)/r$. Hence it is preferable to restrict to metrics with $h'(0) = 0$.  

We report on current work in progress by the second and third authors of this paper with Rubinstein which 
establishes two separate short-time existence results. 
The expected results are as follows: Let $(M^2, g_0)$ be a compact surface with a finite number of conic points, $p_1, \ldots, p_k$, with cone
angles $2\pi \alpha_1, \ldots, 2\pi \alpha_k$. Then there exists a solution $g(t)$ of \eqref{eq:rfs0} on
some small interval $[0,\epsilon)$ for which each $g(t)$ has conic singularities with the same cone angles.
More generally, fix any vector $\gamma = (\gamma_1, \ldots, \gamma_k) \in \RR^k$. Then there
exists a solution $g_\gamma(t)$ defined on some interval $[0,\epsilon)$, where $\epsilon$ depends on $|\gamma|$,
such that each $g_\gamma(t)$ has a $k$-tuple of cone angles $\alpha(t) = (\alpha_1(t), \ldots, \alpha_k(t))$ such that
$\alpha'(0) = \gamma$. 

\medskip

This gives some indication of the extent of nonuniquness. Not only is there an instantaneously complete solution
and a unique one which preserves the cone angles, but in addition, the many others provided by this result 
which change the cone angles. It is not clear whether other solutions exist as well.

We indicate a few ideas that go into the proof. 

The first point is that, as already indicated, it is useful to restrict to conformal factors $\phi(t,r,y)$ of the form
$\phi = \phi_0(t) + \tilde{\phi}(t,r,y)$ where $\tilde{\phi} = \calO(r^2)$. This is possible if all cone angles are less 
than $\pi$; if any cone angle $2\pi \alpha$  lies in $(\pi, 2\pi)$, then we must instead consider
the slightly broader class of functions of the form $\phi = \phi_0(t) + \phi_1(t,y)r^{1/\alpha} + \tilde{\phi}$, 
$\tilde{u} = \calO(r^2)$, where $\phi_1(t,y)$ is a linear combination of $\sin (y/\alpha)$ and $\cos (y/\alpha)$
with coefficients depending on $t$. The functions $\phi_0$ and $\phi_1$ have only finite H\"older regularity in $t$. 
For any function of this form, $\phi(t,r,y) g_0$ has the same cone angles as $g_0$. 

The main technical issue in proving the first short-time existence result is showing that this subclass of 
functions is preserved under propagation by the equation \eqref{eq:rfs}.  (An alternate approach to this 
result (for the flow which keeps cone angles fixed) appears in a paper by Yin \cite{Yin}, who uses an approximation 
argument.  However, his approach does not seem to provide sufficient control of the asymptotics of solutions near the 
cone points, which is very useful to have for many applications, including long-time existence). 

The other point we  discuss now is how the flow can be set up to change the cone angles. For this we
briefly recall how to specify all self-adjoint extensions of the Laplacian $\Delta_g$ for any conic metric $g$
in two dimensions. (Actually, the linearization of the flow equation for $\phi = \log u$ is of the
form $\Delta_g + R$, but since $R$ is bounded, the facts below are the same for this slightly more general operator.)
We begin with a local regularity result: Suppose that $f \in L^2(M^2)$ and $\phi \in L^2(M^2)$
is a solution to $\Delta_{g_0} \phi = f$. Then 
\[
\phi = \phi_0 + \tilde{\phi}_0 \log r + \tilde{\phi}, \qquad \tilde{\phi} = \calO(r^\tau), \ \tau > 1.
\]
(The precise regularity of $\tilde{\phi}$ is not important here.) This solution 
is not unique, and we can specify various solutions by imposing conditions on the coefficients
$\phi_0$ and $\tilde{\phi}_0$. The Friedrichs extension of $\Delta_{g}$ is the operator obtained from imposing
the condition that $\tilde{\phi}_0 = 0$. (In other words, the set of functions $\phi$ with no $\log r$ term and
with $\Delta \phi \in L^2$ comprise a domain on which $\Delta$ is self-adjoint.) However, there are many
other self-adjoint extensions, and these are parametrized by an `angle' $\theta$: We say that
$\phi$ lies in the domain $\calD_\theta$ if $\Delta_{g}\phi \in L^2$ and $\sin \theta \, \phi_0 + \cos \theta 
\, \tilde{\phi}_0 = 0$. (Hence the Friedrichs extension has domain $\calD_0$.)  The possible choices
of solution to $\Delta \phi = f$ when $f$ lies in a H\"older space can be described similarly. 

These are all statements about elliptic theory, but there are analogues for the associated parabolic problem.
Using these to analyze the linearization of the flow, then applying the usual contraction mapping argument,
is what leads to the short-time existence result stated above. Note that the modification of the logarithm of the
conformal factor by a term which includes a (small!) multiple of $\log |z|$ is precisely what changes the cone angles.
By contrast, the instantaneously complete flow produces a solution $u(z,t)$ which has as its leading term a multiple 
of $1/r^2 (\log r)^2$ for any $t > 0$. Equivalently, the logarithm of the conformal factor in that case results in a leading term
of the form $-2 \log r - 2 \log \log r$. 

A more precise description would take us too far afield, so we refer to the forthcoming paper
by Mazzeo, Rubinstein and Sesum for more details, as well as for results concerning the long-time
existence and convergence properties of this flow.


\end{document}